\documentclass[12pt,sloppy]{article}
\usepackage{graphicx}
\usepackage{subfigure}
\usepackage{latexsym}
\usepackage{amsfonts}
\usepackage{amsthm}
\usepackage{amsmath}
\usepackage{caption}
\usepackage{chngpage}
\author{{\bf Maya Mohsin Ahmed }}
\title{{\LARGE {\bf Demystifying  Benjamin Franklin's   other 8-square   }}}
\date{}

\newtheorem{defn}{Definition}[section]

\begin{document}     % this begins the actual document body
\maketitle           % this actually creates the title block

\begin{abstract}
 In   this  article,   we  reveal   how   Benjamin  Franklin
 constructed his  second $8 \times 8$ magic  square. We also
 construct two new $8 \times 8$ Franklin squares.
\end{abstract}

\section{Introduction}

%%%% Figure of Franklin squares %%%%%%%%%%%%%%%%%%%%%%%
\begin{figure}[h]
 \begin{center}
\caption{Squares constructed by Benjamin Franklin.} \label{franklinsquares}
     \includegraphics[scale=0.5]{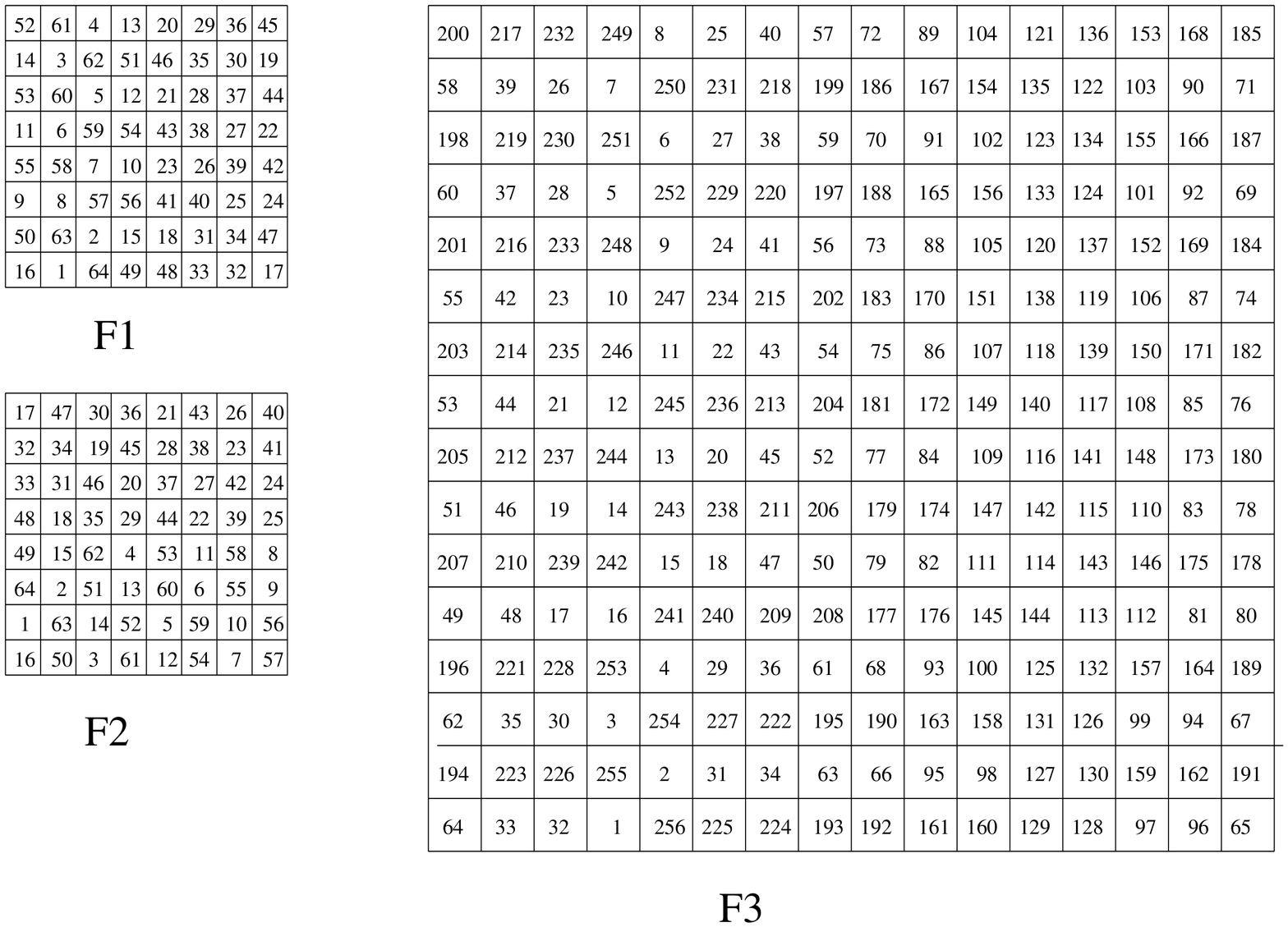}

 \end{center}
 \end{figure}

{\tiny
\captionof{table}{Franklin squares constructed using Hilbert basis \cite{ma}.} \label{n1n2n3table}
%\begin{adjustwidth}{-1in}{.5in}
\[
\begin{array}{ccc}
\begin{array}{|l|l|l|l|l|l|l|l|} \hline
3 & 61 & 16 & 50  & 7 & 57 &  12 &  54 \\  \hline
32 & 34 & 19 & 45 &  28 & 38 & 23 & 41 \\ \hline
33 & 31 &  46 &  20 & 37 & 27 & 42 & 24 \\ \hline
62 & 4 & 49 & 15 & 58 & 8 & 53 & 11 \\ \hline
35 & 29 & 48 & 18 & 39 & 25 & 44 &  22 \\ \hline
64 & 2 & 51 & 13 &  60 & 6 & 55 &  9 \\ \hline
1 & 63 & 14 & 52 & 5 & 59 & 10 &  56 \\ \hline
30 & 36 & 17 & 47 &  26 & 40 & 21 & 43 \\ \hline
\end{array} & 

\begin{array}{|l|l|l|l|l|l|l|l|} \hline
17 & 47 & 30 & 36 & 25 & 39 & 22 & 44 \\ \hline
32 & 34 & 19 & 45 &  24 & 42 & 27 & 37 \\ \hline
33 & 31 & 46 & 20 & 41 &  23 & 38 & 28 \\ \hline
48 & 18 &  35 &  29 & 40 & 26 & 43 & 21 \\ \hline
49 & 15 & 62 & 4 & 57 & 7 & 54 & 12 \\ \hline
64 & 2 & 51 & 13 & 56 & 10 & 59 &  5 \\ \hline
1 & 63 & 14 & 52 & 9 & 55 & 6 & 60 \\ \hline
16 & 50 &  3 &  61 & 8 & 58 & 11 &  53 \\ \hline
\end{array}
 &
\begin{array}{|l|l|l|l|l|l|l|l|} \hline
44 & 61 & 4 & 21 & 12 &  29 & 36 & 53 \\ \hline
22 & 3 & 62 & 43 & 54 & 35 & 30 & 11 \\ \hline
45 & 60 & 5 & 20 &  13 & 28 & 37 &  52 \\ \hline
19 & 6 & 59 & 46 & 51 & 38 & 27 &  14 \\ \hline
47 & 58 & 7 &18 &15 &26 & 39 & 50 \\ \hline
17 & 8 & 57 & 48 & 49 & 40 & 25 & 16 \\ \hline
42 & 63 & 2 & 23 & 10 & 31 & 34 &  55 \\ \hline
24 & 1 & 64 & 41 & 56 & 33 & 32 & 9 \\ \hline
\end{array} \\ \\
\mbox{N1}  & \mbox{N2} & \mbox{N3} 

\end{array}
\]
%\end{adjustwidth}
}

{\tiny
\captionof{table}{New Franklin squares.} \label{brandnewsquaretable}
\[\begin{array}{cc}
\begin{array}{|l|l|l|l|l|l|l|l|} \hline
29 &  35 & 18 & 48 & 21 & 43 & 26 & 40 \\ \hline
20 & 46 & 31 & 33 &  28 & 38 & 23 & 41 \\ \hline
45 & 19 & 34 &  32 & 37 & 27 & 42 & 24 \\ \hline
36 & 30 & 47 & 17 & 44 &  22 & 39 & 25 \\ \hline
61 & 3 & 50 & 16 & 53 &  11 & 58 & 8 \\ \hline
52 & 14 & 63 & 1 & 60 & 6 & 55 & 9 \\ \hline
13 & 51 & 2 & 64 &  5 & 59 & 10 & 56 \\ \hline
4 & 62 & 15 & 49 &  12 & 54 & 7 & 57 \\ \hline
\end{array} &
\begin{array}{|l|l|l|l|l|l|l|l|} \hline
29 & 35 & 18 & 48 & 25 & 39 & 22 & 44 \\ \hline
20 & 46 & 31 & 33 & 24 & 42 & 27 & 37 \\ \hline
45 & 19 & 34 & 32 & 41 & 23 & 38 & 28 \\ \hline 
36 & 30 & 47 & 17 & 40 & 26 & 43 & 21 \\ \hline
61 & 3 & 50 & 16 &  57 & 7 & 54 & 12 \\ \hline
52 & 14 & 63 & 1 & 56 & 10 & 59 & 5 \\ \hline
13 & 51 & 2 & 64 &  9 & 55& 6 & 60 \\ \hline
4 & 62 & 15 & 49 & 8 & 58 & 11 & 53 \\ \hline
\end{array} \\ \\
\mbox{B1} & \mbox{B2}
\end{array}
\]
}

The well-known squares  in Figure \ref{franklinsquares} were
constructed  by   Benjamin  Franklin.  The   square  F2  was
introduced separately  and hence is generally  known as {\em
  the other  8-square}. The entries of the  squares are from
the set $\{1,2, \dots,  n^2\}$, where $n=8$ or $n=16$. Every
integer in this  set occurs in the square  exactly once. For
these squares, the entries of  every row and column add to a
common  sum called  the {\em  magic sum}.  The $8  \times 8$
squares have magic sum 260 and the $16 \times 16$ square has
magic  sum 2056.  In  every  half row  and  half column  the
entries add to  half the magic sum. The  entries of the main
bend diagonals and all the bend diagonals parallel to it add
to the magic sum. In  addition, observe that every $2 \times
2$ sub-square in  F1 and F2 adds to half  the magic sum, and
in F3 adds to one-quarter the magic sum. The property of the
$2  \times 2$  sub-squares adding  to a  common sum  and the
property of bend diagonals adding  to the magic sum are {\em
  continuous  properties}. By  continuous  property we  mean
that if  we imagine  the square is  the surface of  a torus;
i.e. opposite  sides of the square are  glued together, then
the bend  diagonals or the  $2 \times 2$ sub-squares  can be
translated and  still the corresponding sums  hold. In fact,
these squares have  innumerable fascinating properties. See
\cite{ma}, \cite{andrews},  and \cite{pasles} for  a detailed
study of these squares.

From now on, row sum, column sum, or bend diagonal sum, etc.
mean that we are adding the entries of those elements. Based
on the descriptions of  Benjamin Franklin, we define {\em
  Franklin squares} as follows (see \cite{ma2}).

\begin{defn}[Franklin Square]Consider
 an integer, $n=2^r$ such that $r \geq 3$. Let the magic sum
 be  denoted by  $M$ and  $N=n^2+1$.  We define  an {\em  $n
   \times n$  Franklin square} to  be a $n \times  n$ matrix
 with the following properties:
\begin{enumerate}
\item  Every  integer from  the  set  $\{1,2, \dots,  n^2\}$
  occurs  exactly once  in the  square. Consequently,  \[M =
  \frac{n}{2} N.\]
\item All  the the half  rows, half columns add  to one-half
  the magic sum. Consequently,  all the rows and columns add
  to the magic sum.
\item All the bend diagonals add to the magic sum, continuously.

\item All  the $2  \times 2$ sub-squares  add to  $4M/n =2N$,
  continuously.   Consequently,  all   the   $4  \times   4$
  sub-squares add to $8N$, and the four corner numbers
  with the four middle numbers add to $4N$.
\end{enumerate}
\end{defn}

 We call  permutations of the  entries of a  Franklin square
 that  preserve  the  defining  properties of  the  Franklin
 squares  {\em  symmetry operations},  and  two squares  are
 called {\em isomorphic} if we can get one from the other by
 applying symmetry  operations. Rotation and  reflection are
 clearly  symmetry   operations.  See  \cite{ma}   for  more
 symmetry  operations  like   exchanging  specific  rows  or
 columns. In  \cite{ma}, we  showed that F1  and F2  are not
 isomorphic  to each  other. 

Constructing  Franklin  squares  are  demanding and  only  a
handful such  squares are  known to date.  We showed  how to
construct  F1, F2,  and F3  using methods  from  Algebra and
Combinatorics  in \cite{ma}.  In the  same article,  for the
first  time  since  Benjamin  Franklin, we  constructed  new
Franklin   squares  N1,   N2,   and  N3,   given  in   Table
\ref{n1n2n3table}.  We proved  that these  squares  were not
isomorphic to  each other nor to  F1 or F2.  In other words,
these   squares  were  really   new.  Those   methods  being
computationally  challenging  are  not suitable  for  higher
order  Franklin squares.  Moreover,  the constructions  used
computers  and  hence  lacked  the  intrigue  of  Franklin's
constructions. In \cite{ma2},  we followed Benjamin Franklin
closely,  and  used  elementary  techniques to  construct  a
Franklin square of every order. With these techniques we are
able  to  construct  F1  and  F3, but  not  F2.  In  Section
\ref{constructsection}, we modify  the methods in \cite{ma2}
to  construct F2.  In Section  \ref{newconstructsection}, we
create  new  Franklin squares  B1  and  B2,  given in  Table
\ref{brandnewsquaretable},  but this  time  using elementary
techniques,  in keeping  with  the true  spirit of  Benjamin
Franklin.

\section{Franklin's construction of his other $8$-square.} \label{constructsection}
In  this section,  we  show how  to  construct the  Franklin
square F2.

\begin{figure}[h]
 \begin{center}
     \includegraphics[scale=0.45]{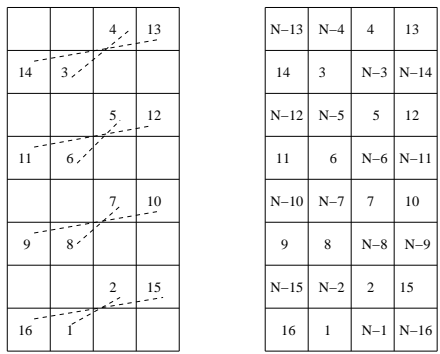}
\caption{Construction of the left half of the Franklin square F1.} \label{f1leftconstruct}
 \end{center}
 \end{figure}
%%%%%%%%%%%%%%%%%%%%%%%%%%%%%%%%%%%%%%%%%%%%%%%%

\begin{figure}[h]
 \begin{center}
     \includegraphics[scale=0.45]{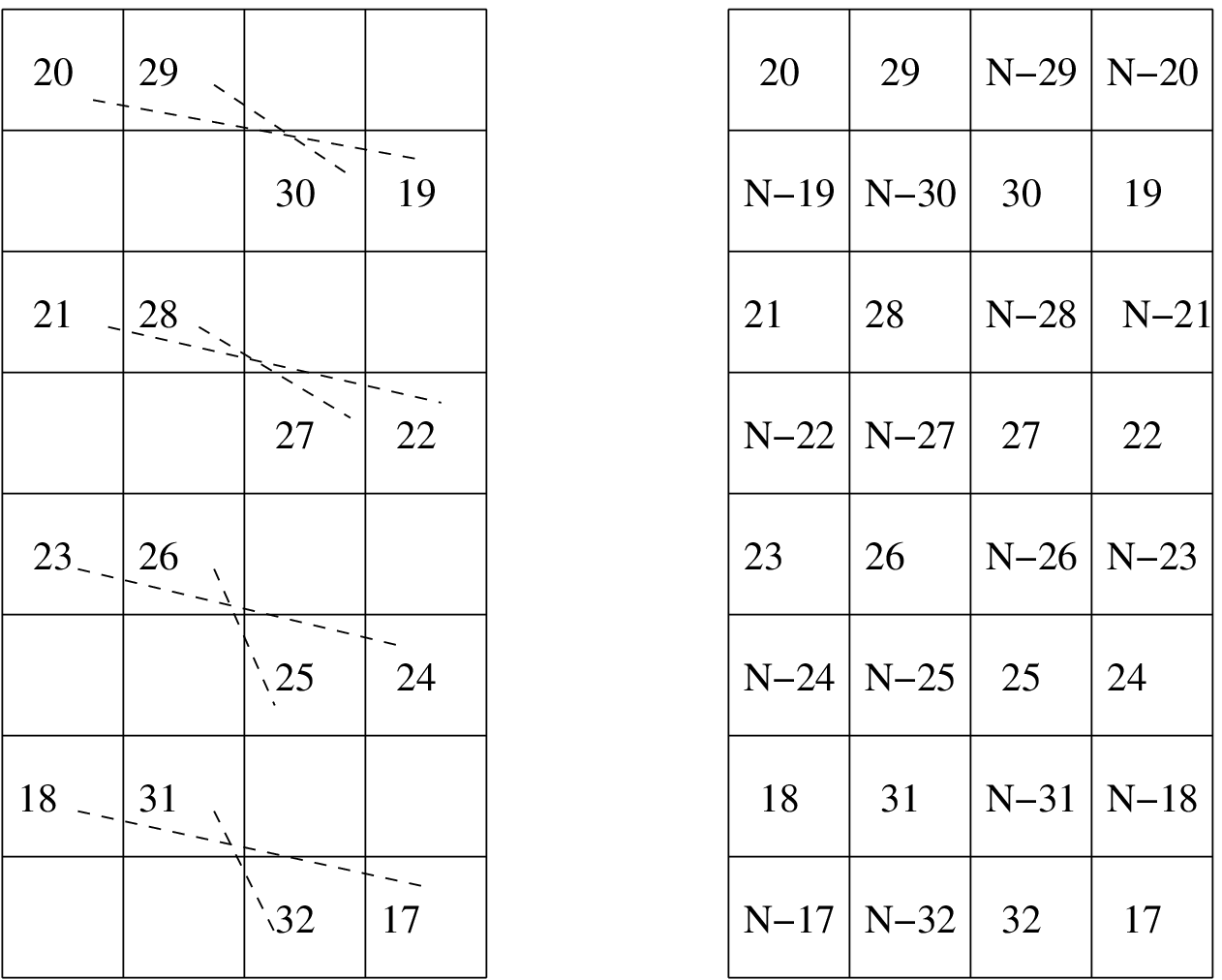}
\caption{Construction of the right half of the Franklin square F1.} \label{f1rightconstruct}
 \end{center}
 \end{figure}
%%%%%%%%%%%%%%%%%%%%%%%%%%%%%%%%%%%%%%%%%%%%%%%%

We  present a  quick  recap  of our  construction  of F1  in
\cite{ma2}. Note that $N=65$ throughout this article. The $8
\times 8$ square was divided in to two halves: the left half
consisting  of  first  four   columns  and  the  right  half
consisting  of  the last  four  columns.  The  left half  is
constructed by  first placing the numbers $1,  2, \dots, 16$
in    a    specific    pattern    as   shown    in    Figure
\ref{f1leftconstruct}. Then,  the number $N-i$  is placed in
the same row that contains $i$, where $i=1, \dots, 16$, such
that $i$ and $N-i$ are always equidistant from the middle of
the left half.  See Figure \ref{f1leftconstruct}. Similarly,
the second half is constructed by first placing numbers $17,
18, \dots,  32$ in a specific pattern.  Finally, the numbers
$N-i$ is  placed in  the same row  that contains  $i$, where
$i=17,  \dots,  32$, such  that  $i$  and  $N-i$ are  always
equidistant  from  the  middle   of  the  right  half.  This
procedure is explained in Figure \ref{f1rightconstruct}. The
algorithm for generating the  pattern of $1,2, \dots, 32$ in
F1  is  given  in   \cite{ma2}.  We  also  generalized  this
procedure to construct a Franklin square for any given order
in \cite{ma2}. From now  on, throughout the article, when we
say {\em  pattern of a square},  we mean the  pattern of the
numbers $1,2, \dots, 32$ in the square.

\begin{figure}[h]
 \begin{center}
     \includegraphics[scale=0.45]{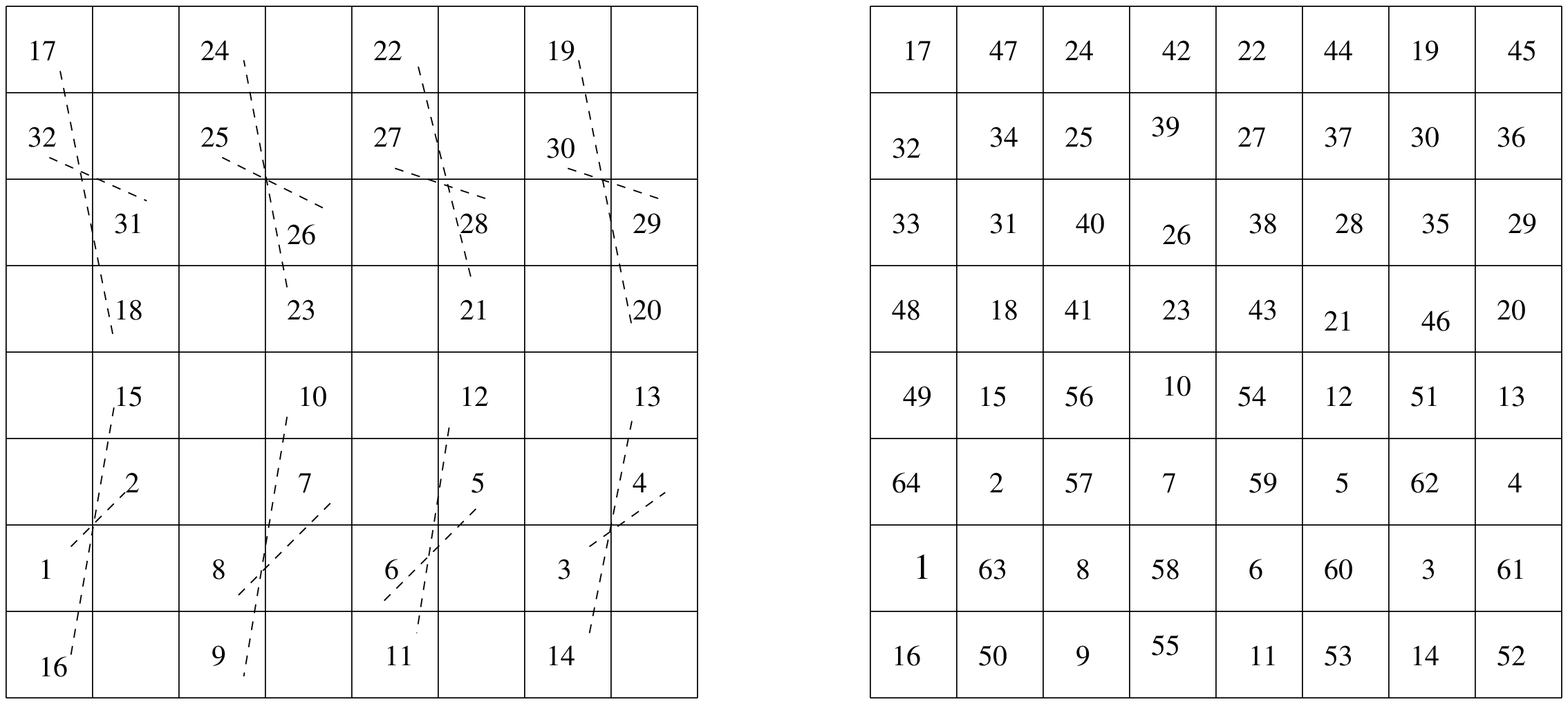}
\caption{Construction of Franklin square S.} 
\label{columnsconstruct}
 \end{center}
 \end{figure}
%%%%%%%%%%%%%%%%%%%%%%%%%%%%%%%%%%%%%%%%%%%%%%%%

\begin{figure}[h]
 \begin{center}
     \includegraphics[scale=0.45]{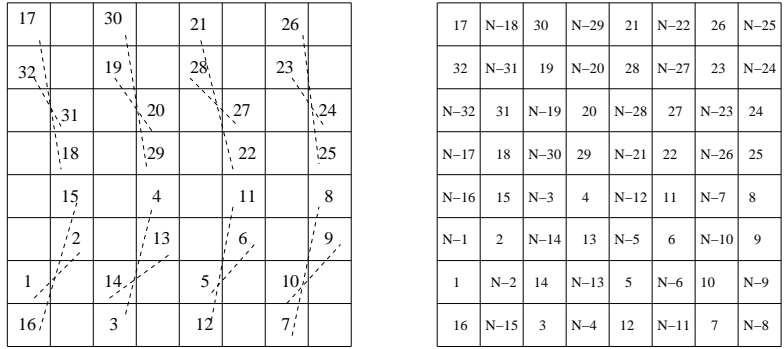} 
\caption{ Constructing Franklin square F2.}
\label{f2constructfig}
 \end{center}
 \end{figure}

In our  search for  new Franklin squares,  it is  natural to
swap rows for  columns in the above method.  Thus, we divide
the $8  \times 8$  square in  to two halves  : the  top half
consisting  of  the first  four  rows  and  the bottom  half
consisting of the  last four rows. We place  the numbers $1,
\dots, 16$ in  the bottom half, and the  numbers $17, \dots,
32$  in  the  top  half  in  the  pattern  shown  in  Figure
\ref{columnsconstruct}. This time we place the numbers $N-i$
in the same  columns as $i$, equidistant from  the middle of
each part. The Franklin square  S we get from this procedure
is isomorphic  to F1: check  that we can derive  this square
from F1 by rotating it  and then reflecting it. The square S
is not interesting because it is essentially the same as F1.
But the pattern of entering  the numbers $1,2, \dots, 32$ in
S  is algorithmic  since it  is just  a row  version  of the
pattern of F1. That is,  instead of filling two columns at a
time, we fill two rows at a time. We modify the pattern of S
to get a  pattern for F2. After that, we  place $N-i$ in the
same column as  $i$, for $i=1, \dots 32$,  equidistant from the
center    of    each   half,    as    usual.   See    Figure
\ref{f2constructfig} for the construction of F2.

%%%%%%%%%%%%%%%%%%%%%%%%%%%%%%%%%%%%%%%%%%%%%%%%

%%%%%%%%%%%%%%%%%%%%%%%%%%%%%%%%%%%%%%%%%%%%%%%%

The construction of F2 involved moving columns of in the $1,
\dots,  32$  pattern of  the  square  S  and then  switching
elements  along the  diagonals. This  pattern is  not easily
generalized to  higher orders. But  in the next  Section, we
modify the patterns of known squares to create new squares.

\section{Constructing new Franklin squares.} \label{newconstructsection}

In this section we construct the new Franklin squares B1 and B2 (see
Table \ref{brandnewsquaretable}).

\begin{figure}[h]
 \begin{center}
     \includegraphics[scale=0.45]{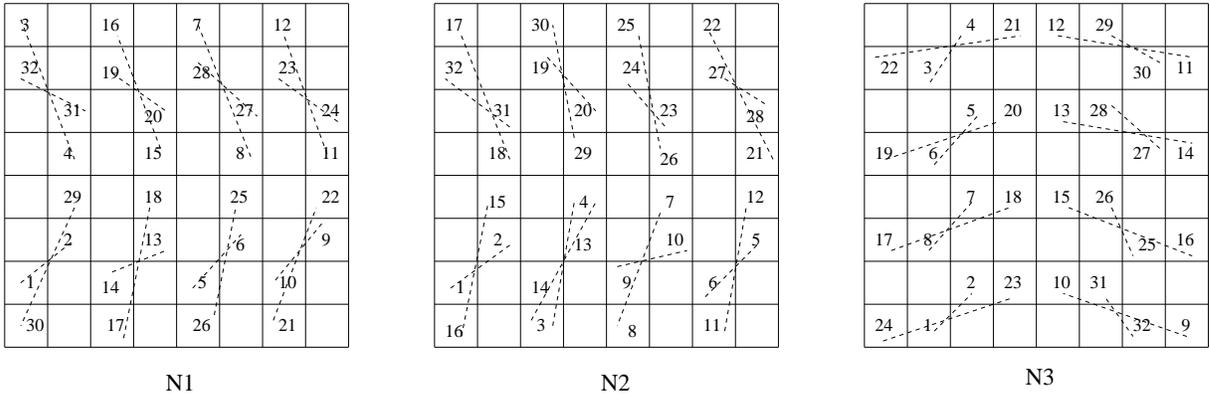}
\caption{The patterns of  Franklin squares
  N1, N2, and N3.}
\label{nipatternsfig}
 \end{center}
 \end{figure}
%%%%%%%%%%%%%%%%%%%%%%%%%%%%%%%%%%%%%%%%%%%%%%%%

To begin with, we look at the patterns of the numbers $1, 2,
\dots, 32$ appearing in the Franklin squares N1, N2, and N3.
These are shown in Figure \ref{nipatternsfig}. We found that
even  in these  squares that  were constructed  using Hilbert
bases (see  \cite{ma}), the strategy of  finding the pattern
of the  Franklin square, and  then placing $i$ and  $N-i$ in
the same row or column, as the case may be, equidistant from
the center  of relevant half of the  square, worked. Observe
that the patterns of N1  and N2 are derived by modifying the
pattern of F2. So in  their constructions, we will place $i$
and $N-i$ in the same  columns. On the other hand, since, N3
is  a permutation of  the pattern  of F1,  we place  $i$ and
$N-i$  in  the same  rows  while  constructing N3.  Benjamin
Franklin's  patterns  always  restrict  the entries  $1,  2,
\dots, 16$ to one half of the square. Observe that N2 is the
only square with this property.

\begin{figure}[h]
 \begin{center}
     \includegraphics[scale=0.45]{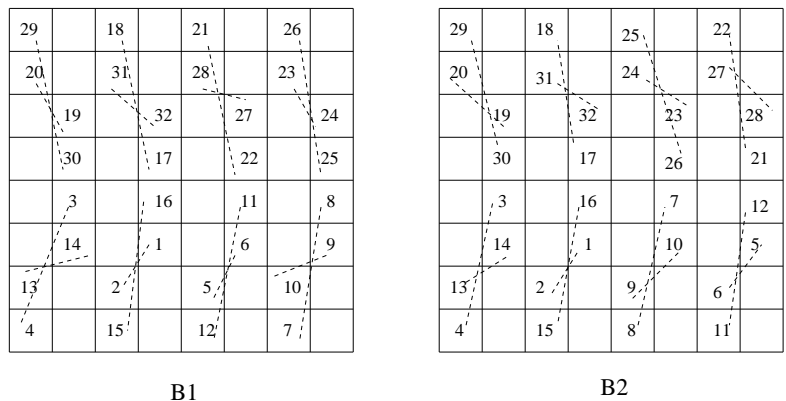}
\caption{The patterns of the Franklin squares B1 and B2.}
\label{b1b2pattern} 
\end{center}
 \end{figure}

Guided by these observations,  we modified the pattern of F2
to  build the  new Franklin  squares B1  and B2.  See Figure
\ref{b1b2pattern} for  the patterns of B1 and  B2. Since the
patterns were  derived from F2,  to construct them,  we will
place $i$ and  $N-i$ in the same columns. B1  and B2 are not
isomorphic to each other  or any other known squares because
the columns were made different by construction.
%%%%%%%%%%%%%%%%%%%%%%%%%%%%%%%%%%%%%%%%%%%%%%%%

Clearly, there  are more Franklin squares. To  find them, we
need  to  find more  patterns  that  will  yield a  Franklin
square. So  it is time to  ask again ``How  many squares are
there Mr. Franklin?''

\end{document}